\documentclass[letterpaper,12pt]{article}

\usepackage[colorlinks, linkcolor=blue]{hyperref}

\usepackage{amsfonts,amsmath,verbatim,graphicx, subfigure,longtable}
\usepackage{amssymb,amsthm,amscd}
\newtheorem{thm}{Theorem}[section]

\newtheorem{prop}[thm]{Proposition}

\newtheorem{definition}[thm]{Definition}

\newtheorem{example}[thm]{Example}
\newtheorem{remark}[thm]{Remark}

\newenvironment{pf}{\noindent\emph{Proof \,}}{\mbox{}\qed}

\newcommand{\toL}{\,{\buildrel \mathcal{L} \over \longrightarrow}\,}
\newcommand{\equalL}{\,{\buildrel \mathcal{L} \over =}\,}

\numberwithin{equation}{section}
\def\R{{\mathbb R}}     % Real numbers
\def\E{{\mathbb E}}     % Expectation
\def\X{{\mathfrak X}}   % Empirical process
\def\Y{{\mathcal Y}}   % Fluctuation process

\def\L{{\mathcal L}}   % Generator of \eta

\def\eps{\varepsilon}
\def\<{{\langle}}
\def\>{{\rangle}}
\def\1{{\bf 1}}         % Indicator

\renewcommand{\bar}{\overline}

\begin{document}
\allowdisplaybreaks

\title{\Large \bf
Correlation function methods for a system of  annihilating Brownian particles
\thanks{Mathematics Department, Indiana University. Email: \texttt{waifan@iu.edu}. Research partially supported by NSF Award DMS-1804492.}
%\thanks{Research partially supported by NSF Award DMS-1804492.}
 }

\author{{\bf Wai-Tong (Louis) Fan}}
\date{\today}
\maketitle

	\begin{abstract} 
In this expository note we highlight the correlation function method as a unified approach in proving both hydrodynamic limits and fluctuation limits for reaction diffusion particle systems. 
%We review the key steps in establishing the functional law of large numbers for a simple stochastic particle system. 
For simplicity we focus on the case when the hydrodynamic limit is $\partial_t u=\frac{1}{2}\Delta u -u^2$, one of the simplest nonlinear  reaction-diffusion equations.  The outline of the proof follows from Chapter 4 of De Masi and Presutti \cite{MP91} but to simplify the presentation, we consider reflected Brownian motion instead of reflected random walks.  We also briefly mention the key ideas in proving the fluctuation result.

%Key ideas about using correlation functions to prove fluctuation results are also mentioned. 
% 	The aim of this handout is to highlight the method of correlation functions as a unified approach in proving both hydrodynamic limit and fluctuation limit for stochastic particle systems.
	\end{abstract}

\bigskip
\noindent {\bf AMS 2000 Mathematics Subject Classification}: Primary
60F17, 60K35; Secondary  92D15

\bigskip\noindent
{\bf Keywords and phrases}: correlation functions, hydrodynamic limit, interacting diffusion, reflected diffusion, annihilation, non-linear partial differential equation

\section{Introduction}

It is known that partial differential equations (PDE) and stochastic partial differential equations (SPDE) can arise respectively as hydrodynamic limits and fluctuation limits of interacting particle systems. These results,  often formulated as 
functional law of large numbers (LLN) and functional central limit theorems (CLT),
 are very desirable  for various reasons. First, they are universal in the sense that the limits are robust against fine details of the underlying particle systems. This advantage is typically not carried over to large deviation results such as \cite{BFW}. Second, the hydrodynamic limit $u$ often gives the asymptotic behavior of the joint law of a fixed number of particles. For example, for exchangeable systems,  the LLN is equivalent to the  propagation of chaos (POC) \cite{Szn91} and the latter says that for any fixed $k$, the joint law of $k$ particles is given by the product $u^{\otimes k}$,  as the number of particles tends to infinity.

One of the most extensively studied stochastic particle systems are  those associated with reaction-diffusion equations of the form
\begin{equation}\label{E:RDeqt}
  \partial_t u(t,x)=\frac{1}{2}\Delta u(t,x) + R(u(t,x)),
\end{equation}
where $\Delta$ is the standard Laplacian representing diffusion of a population and $R(u)$ is a function in $u$, called the reaction term, representing a certain interaction in the population. 

An elegant example of such a particle system is studied in P. Dittrich \cite{pD88a}. One starts with $N$ particles on the unit interval $[0,1]$ which perform reflecting Brownian motions (RBMs) and specifies that, for $2\leq k\leq n$, any $k$-tuple of particles with pairwise distances $\eps=1/N$, say $(x^{i_1},\cdots,\,x^{i_k})$, disappears simultaneously with intensity
$$c_k(k-1)! \eps^{k-1}\,\int_{[0,1]}p(\eps^2,x^{i_1},y)\cdots p(\eps^2,x^{i_k},y)\,dy,$$
where $p(t,x,y)$ is the transition density of a RBM on $[0,1]$ and $c_k\in [0,\infty)$.
It is shown that, as $N\to\infty$, the hydrodynamic limit  is described by \eqref{E:RDeqt} with Neumann boundary condition and 
$$R(u)=-\sum_{k=2}^n c_k u^k.$$ 

Fluctuation results are briefly discussed in Section \ref{S:fluctuation}. The focus for now is to prove the hydrodynamic limit result for the special case $R(u)=-u^2$. The proof follows Chapter 4 of De Masi and Presutti \cite{MP91}. However, to simplify the exposition we use RBMs rather than random walks as microscopic dynamics. 
The principles of a correlation function technique (which involve the BBKGY hierarchy\footnote{BBGKY stands for N. N. Bogoliubov, Max Born, H. S. Green, J. G. Kirkwood, and J. Yvon, who derived this type of  hierarchy of equations in the 1930s and 1940s in a series of papers.}) are highlighted.

\section{The model}\label{S:model}
Intuitively speaking, the aforementioned process consists of $N$ independent RBMs on $[0,1]$ and any pair of them, say $(x,y)$, disappears with intensity $(1/N)\,p(2/N^2,x,y)$. The Gaussian estimates for $p(t,x,y)$ implies that the interaction distance is of order $1/N$, and the intensity of disappearance for such a pair is of order 1. This scaling is reasonable since the number of pairs is of order $N$ (imagine all particles are spread out evenly on the interval).

Precisely, we construct a family (indexed by $N$) of Markov processes by specifying their infinitesimal generators as follow. Let $S_m$ be the space of unordered $m$-tuples of elements in $[0,1]$ for $m\geq 1$ and $S_0$ be an abstract point representing an absorbing state (when all particles die out). Let $\mathbb{X}^N=(\mathbb{X}^N_t)_{t\geq 0}$ be the continuous time Markov process with state space $\mathbb{S}:= \cup_{m=0}^{\infty}S_m $ and with generator $\L_N$ defined by
\begin{equation}\label{E:GeneratorL_N}
\L_Nf(\bar{x}):= \frac{1}{2}\sum_{i=1}^{m}\frac{\partial^2}{(\partial x^i)^2}f(\bar{x})+ \frac{1}{2N}\sum_{i\neq j}[f(\bar{x}_{ij})-f(\bar{x})]\,p(2/N^2,x,y),\quad\bar{x}\in S_m,\,m\geq 1,
\end{equation}
where $\bar{x}_{ij}$ is the element of $S_{m-2}$ obtained from $\bar{x}=(x^1,\cdots,\,x^m)\in S_m$ by deleting $x^i$ and $x^j$. When $\bar{x}\in S_0$, we define $\L_Nf(\bar{x}):=0$. %See Remark \ref{Welldefine} for well-defineness of $\mathbb{X}^N$.

Such a process $\mathbb{X}^N$ is well-defined:
the domain of $\L_N$ contains the class of functions $f:\,\mathbb{S}\to\R$ whose restrictions to $S_m$ belong to $C_2([0,1]^m)$ for all $m\geq 1$, where $C_2([0,1]^m)$ denotes the space of twice continuous differentiable functions on $[0,1]^m$ whose normal derivatives vanish on the boundary of $[0,1]^m$. It is then routine to check that all hypothesis of the Hill-Yosida theorem (see, for instance, Chapter 1 of \cite{EK86}) are satisfied. Hence  $\L_N$ determines a unique Markov process in distribution.

\begin{remark}\rm
	In this note we consider only ``soft annihilation" in which annihilation occurs with a probability less than 1. We could have instead specify ``hard annihilation" in which annihilation occurs with probability one when two particles are within an interaction distance $\delta_N$. Results of Sznitman \cite{Sznitman87} suggest that the same LLN limit \eqref{E:RDeqt2} can be obtained if	$\delta_N$ is of order $N^{-1/(d-2)}$ when $d\geq 3$ and of order $e^{-N}$ when $d=2$.
\end{remark}

\section{Functional law of large numbers}

For each $N$ and $t\geq 0$, we have either $\mathbb{X}^N_t=(x^1_t,\cdots,\,x^{m(t)}_t)$ for some positive integer $m(t)=m_N(t)$ (the number of particles alive at time $t$) or $\mathbb{X}^N_t\in S_0$. The \emph{normalized empirical distribution} of the particles alive is
\begin{equation}\label{Def:EmpiricalRD}
    \X^N_t(dz)  := \dfrac{1}{N}\sum_{i=1}^{m(t)}\1_{x^{i}_t}(dz)  \quad \text{if }m(t)\geq 1
\end{equation}
and is defined as the zero measure if $\mathbb{X}^N_t\in S_0$. Note that $\X^N_t$ is a random measure on $\bar{D}=[0,1]$. Moreover, $\X^N= (\X^N_t)_{t\geq 0}$ is a strong Markov process in $M_+(\bar{D})$, the space of finite non-negative Borel measures on $\bar{D}$ equipped with weak topology, and $\X^N$ has sample paths in the Skorokhod space $D([0,\,\infty),\,M_{+}(\bar{D}))$.

In what follows, $\toL$ denotes convergence in probability law, $\equalL$ denotes equal in probability law. For a topological space $E$, we denote by $C(E)$ the space of continuous functions on $E$.

\begin{thm}\label{T:LLN_RD}\textbf{(Functional Law of Large Numbers)}
    Suppose $\{\X^N_0\} \toL u_0(x)\,dx$ in $M_+(\bar{D})$, where $u_0\in C(\bar{D})$. Then
    $$\X^N_t(dx) \,\toL\, u(t,x)\,dx \quad \text{in  }D([0,\,\infty),\,M_{+}(\bar{D})),$$
    where $u\in C([0,\infty)\times\bar{D})$ is the solution 
    to the reaction-diffusion equation
    \begin{equation}\label{E:RDeqt2}
     \partial_t u(t,x)=\frac{1}{2}\Delta u(t,x) - u^2(t,x)
    \end{equation}
    with Neumann boundary condition and initial condition $u(0,\cdot)=u_0$.
\end{thm}

\medskip

In Theorem \ref{T:LLN_RD}, $u$ is the unique element in $C([0,\infty)\times\bar{D})$ which satisfies the integral equation
\begin{equation}\label{Mildsol}
u(t,x)= P_tu_0(x)-\int_0^{t}P_{t-s}(u(s,\cdot))(x)\,ds.
\end{equation}
The fact that such $u$ exists can be checked by a fixed point argument. 
$u$ is also called a ``probabilistic solution" (see \cite{zqCwtF14b}) since
it satisfies $u(t,x)= \E^{x}\left[u_0(X_t) \exp \left(-\int_0^tu(t-s,X_s)\,ds \right)\right]$, where $X_t$ is the RBM on $\bar{D}$.

The key technique used in the proof is in the next section.

\section{Propagation of chaos}

\begin{definition}\label{Def:CorrelationFcn}
    Fix $N\in \mathbb{N}$ and consider the process $\mathbb{X}^N_t=(x^1_t,\cdots,\,x^{m(t)}_t)$ constructed above. For $k\geq 1$ and $t\geq 0$, the $k$\textbf{-correlation function at time} $t$, $F^{(k)}_t=F^{N,(k)}_t$, is define to be the function (up to Lebesque a.e.) satisfying
    \begin{equation*}
    \int_{D^k}\Phi(\vec{x})\,F^{(k)}_t(\vec{x})\,d\vec{x}=
    \E\bigg[
        \dfrac{1}{N^{(k)}}\sum_{\substack{i_1,\cdots, i_n\\ \text{distinct}}}^{m(t)}\,\Phi(x^{i_1}_t,\cdots, \,x^{i_k}_t)
    \bigg] \quad\text{for all }\Phi\in C(\bar{D}^k),
    \end{equation*}
    where $N^{(k)}:= N(N-1)\cdots (N-k+1)$ is the number of permutations of $k$ objects chosen from $N$ objects.
\end{definition}

Intuitively, if we randomly pick $k$ living particles in $D$ at time $t$, then $F^{(k)}_t(\vec{x})$ is the probability joint density function for their positions. Note that $F^{(k)}_t$ is defined for almost all $\vec{x}\in D^k$, and that it depends on both $N$ and the initial configurations. We will see that $F^{(k)}_t\in C(\bar{D}^k)$ for $t>0$. We can also replace $N^{(k)}$ by $N^k$ because the behavior of $F^{k}$ when $N\to\infty$ is our concern and $N^{(k)}/N^k \nearrow 1$  as $N\to\infty$.

\begin{example}\label{Eg:CorrelationFcn}
Let $\<\X^N_t,\phi\> = \frac{1}{N}\sum_{i=1}^{m(t)}\phi(x^{i}_t)$ be the integral of a test function (an observable) $\phi$ with respect to  $\X^N_t$. Then 
\begin{eqnarray*}
\E[\<\phi,\X^{N}_t\>]&=&\int_{D}\phi(x)\,F^{(1)}_{t}(x)\,dx \quad\text{and}\\
\E[\<\phi,\X^{N}_t\>^2]&=&\frac{1}{N}\int_D\phi^2(x)\,F^{(1)}_t(x)\,dx + \frac{N-1}{N}\int_{D^2}\phi(x_1)\phi(x_2)\,F^{(2)}_t(x_1,x_2)\,dx_1\,dx_2.
\end{eqnarray*}
\end{example}

\medskip
\begin{thm}\label{T:POC_RD}\textbf{(Propagation of Chaos)}
Suppose $\{\X^N_0\} \toL u_0(x)\,dx$ in $M_+(\bar{D})$, where $u_0\in C(\bar{D})$. Then for any $k\geq 1$, we have
\begin{equation}\label{E:POC_RD}
   \lim_{N\to\infty} \sup_{\substack{t\in[0,T]\\(x_1,\cdots,x_k)\in\bar{D}^k}}\Big|F^{N,(k)}_t(x_1,\cdots,x_k)\,-\,\prod_{i=1}^ku(t,x_i)\Big|\,=\,0,
\end{equation}
$u$ is the solution 
to the reaction-diffusion equation \eqref{E:RDeqt2} in Theorem \ref{T:LLN_RD}.
\end{thm}

\begin{pf}(Sketch)
\textbf{Step 1: BBGKY hierarchy for $F^{N,(k)}$. }
Since the interactions of our process only involves annihilations, it is immediate that $F^{(k)}_t\leq P^{(k)}F^{(k)}_0$, where $P^{(k)}_t$ is the semigroup for the RBM on $[0,1]^k$. Applying Dynkin's formula to the functional
$$(s,\,\X^{N}_s) \mapsto
\dfrac{1}{N^{(k)}}\sum_{\substack{i_1,\cdots, i_k\\ \text{distinct}}}^{\sharp_s}\,P^{(k)}_{t-s}\Phi(x^{i_1}_s,\cdots, x^{i_k}_s)
\;\;,\,s\in[0,t]$$
yields, via the formula \eqref{E:GeneratorL_N} for $\L_N$, the system of equations
\begin{equation}\label{E:BBGKY_F}
F^{(k)}_t=P^{(k)}_t F^{(k)}_0 - \int_0^t P^{(k)}_{t-s}\left( RF^{(k+1)}_s + \frac{QF^{(k)}_s}{N} \right)\;ds,
\end{equation}
where $R$ and $Q$ are operators defined by
\begin{eqnarray*}
RF^{k+1} (x_1,\cdots,x_k)&:=& \sum_{i=1}^k\int_DF^{(k+1)}(x_1,\cdots,x_k,x_{k+1})\,p(2/N^2,x_k,x_{k+1})\,dx_{k+1}\quad\text{and}\\
QF^{(k)}(x_1,\cdots,x_k)&:=& \sum_{i<j}^k F^{(k)}p(2/N^2,x_i,x_j).
\end{eqnarray*}
The system of equations \eqref{E:BBGKY_F} is called the {\it BBGKY-hierarchy} for the correlation functions $F^{N,(k)}$. It is a finite system with exactly $N$ equations, since $F^{N,(N+i)}$ is a zero function for $i\geq 1$.

\textbf{Step 2: Compactness of $\{F^{N,(k)}\}$. }
Using basic properties of the transition kernel $p(t,x,y)$, we can check that
for any $k\geq 1$, the family of functions $\{F^{N,(k)}\}_{N\geq 1}$ is uniformly bounded and equi-continuous on $\bar{D}^k \times [0,\infty)$.
From the above compactness result, it follows that for any sequence $N' \to \infty$ there is a subsequence $N''$ along which $F^{N'',(k)}$ converges, for every $k\geq 1$, uniformly on $\bar{D}^k \times [0,T]$ to some $\gamma^{(k)}\in C(\bar{D}^k \times [0,T])$.

\textbf{Step 3: Limiting hierarchy. }
It can be justified, by passing to the limit $N\to\infty$ for \eqref{E:BBGKY_F} and using basic properties of the heat kernel $p(t,x,y)$, that $\{\gamma^{(k)}\}_{k\geq 1}$ satisfies the limiting infinite hierarchy
\begin{equation}\label{E:BBGKY_infinite}
\gamma^{(k)}_t(\vec{x})= P^{(k)}_t \gamma^{(k)}_0(\vec{x}) - \sum_{i=1}^k\int_0^t P^{(k)}_{t-s}\Big(\gamma^{k+1}_s(z_1,\cdots,z_k,\,z_i)\Big)(\vec{x})\;ds,
\end{equation}
where $P^{(k)}_{t-s}$ acts on the $\vec{z}$ variables. Moreover, it is easy to check that $\prod_{i=1}^k u(t,x_i)$ also satisfies \eqref{E:BBGKY_infinite}. 

\textbf{Step 4: Uniqueness of limiting hierarchy. }
Finally, we must check that the \emph{infinite} limiting hierarchy \eqref{E:BBGKY_infinite} cannot have two distinct solutions. This follows from an easy Gronwall-type argument, using the uniform norm.  The idea is that the difference of two solutions can be bounded above by a sum of $M!$ iterated integrals each of which is bounded by $(Ct)^M/M!$ where $M$ is the number of iterations. So we have uniqueness for small time. Using the semigroup property, we can extend uniqueness to any finite time horizon.
From Step 3 and Step 4, we have
 $$\gamma_t^{(k)}(\vec{x})=\prod_{i=1}^k u(t,x_i)$$
for all $\vec{x}=(x_1,\cdots,x_k)\in \bar{D}^k$, $t\in[0,T]$ and $k\geq 1$. The proof is complete.
\end{pf}

\begin{remark}\label{Uniqueh}\rm
	%Proving uniqueness of \eqref{E:BBGKY_infinite} is easy. However, u
	Uniqueness of infinite limiting hierarchy is usually challenging to obtain. See \cite{lEbShtY07, zqCwtF13a} which  require choosing suitable norms for the correlation functions and manipulations of  the Feynman diagrams or infinite trees. 
\end{remark}

\section{Proof of functional LLN}

The proof of Theorem \ref{T:LLN_RD} now follows from the $C$-tightness of $\{\X^N\}$. This is because by Theorem \ref{T:POC_RD}, the first two moments of $\<\phi,\,\X^{\infty}_t\>$ are identified for all $t\geq 0$ and $\phi\in C(\bar{D})$, where $\X^{\infty}$ is an arbitrary subsequential limit of $\{\X^N\}$. Precisely, we have the following two propositions.

\begin{prop}\label{prop:Tight_RD}
	For all $T>0$, the sequence $\{\X^N\}$ is tight in $D([0,\,T],\,M_{+}(\bar{D}))$. Moreover, any subsequential limit has continuous path almost surely.
\end{prop}

\begin{pf}
	(Sketch)
	The compact containment condition in \cite{EK86} obviously holds. 
	%As pointed out in Remark \ref{Rk:domain}, 
	The domain of the Feller generator $Dom(\frac{1}{2}\Delta)$ is dense in $C(\bar{D})$. Hence it suffice to show the one-dimensional processes $\{\<\X^N,\,\phi\>\}$ is tight in $D([0,\,T],\,\R)$ for all $\phi\in Dom(\frac{1}{2}\Delta)$.
	
	The key is to write down the martingale representation of $\<\X^N,\,\phi\>$. From \eqref{E:GeneratorL_N}, we have for $\phi\in Dom(\frac{1}{2}\Delta)$,
	\begin{equation}\label{E:Mtg_phi}
	\<\X^N_t,\phi\> = \<\X^N_0,\phi\> + \int_0^t \<\X^N_s,\,\frac{1}{2}\Delta\phi\> - \frac{1}{N^2}\sum_{i\neq j}^{m(s)}p(2/N^2,x^i_s,x^j_s)\,\phi(x^i_s) + M^{\phi}_N(t)\;ds,
	\end{equation}
	where $M^{\phi}_N(t)$ is a martingale with quadratic variation
	\begin{equation}\label{E:Mtg_phi_quad}
	\<M^{\phi}_N\>_t=
	\int_0^t \<\X^N_s,\,|\nabla\phi|^2\> + \frac{2}{N^2}\sum_{i\neq j}^{m(s)}p(2/N^2,x^i_s,x^j_s)\,\Big(\frac{\phi(x^i_s)+\phi(x^j_s)}{2}\Big)^2 \;ds.
	\end{equation}
	We can then check tightness of $\{\<\X^N,\,\phi\>\}$ in $D([0,\,T],\,\R)$ by applying Prohorov's Theorem, using standard estimates of the heat kernel $p(t,x,y)$.
\end{pf}

\begin{prop}\label{prop:Mean_Var_RD}
	For all $\phi\in C(\bar{D})$ and $t\geq 0$, we have
	\begin{eqnarray}
	\E^{\infty}[\<v(t),\,\phi\>]&=&\<u(t),\,\phi\>\quad \text{and}
	\label{E:Mean}\\
	\E^{\infty}[\<v(t),\,\phi\>^2]&=&\<u(t),\,\phi\>^2,
	\label{E:Var}
	\end{eqnarray}
	where $\E^{\infty}$ is the law of an arbitrary subsequential limit $\X^{\infty}$ of  $\{\X^N\}$, and $v(t,x)$ is the density of $\X^{\infty}$, w.r.t. Lebesque measure.
\end{prop}

Equations \eqref{E:Mean} and \eqref{E:Var} follow immediately from Theorem \ref{T:POC_RD} and Example \ref{Eg:CorrelationFcn}.
The proof of Theorem \ref{T:LLN_RD} is complete.

\begin{remark}\rm
	Equation \eqref{E:Mtg_phi} also motivates our choice of the interaction intensity. In fact, it can be any non-negative continuous symmetric function of the form $(1/N)r_N(x,y)$ with $\int_Dr_N(x,y)dy=1$ for all $x$ and $r_N(x,y)\leq C\,p(2/N^2,x,y)$ for all $x,y$ (where $C$ is a positive constant). The choice $(1/N)\,p(2/N^2,x,y)$ is for simplicity. 
	For example, the nonlinear term is
	\begin{eqnarray*}
		&&\frac{1}{N^2}\sum_{i\neq j}^{m(s)}p(2/N^2,x^i_s,x^j_s)\,\phi(x^i_s)\\
	&=&\frac{1}{N^2}\sum_{i, j=1}^{m(s)}p(2/N^2,x^i_s,x^j_s)\,\phi(x^i_s)-\frac{1}{N^2}\sum_{i=1}^{m(s)}p(2/N^2,x^i_s,x^i_s)\,\phi(x^i_s) \\
		&=& \<p(2/N^2,z,w)\phi(z),\,\X^N_s(dw)\otimes \X^N_s(dz)\> - \frac{1}{N}\<p(2/N^2,x,x),\,\X^N_s(dx)\>
	\end{eqnarray*}
	which formally tends to the desired $\int_D \phi(z)\,u^2(s,z)\,dz$. 
	%Justifying this convergence is nontrivial.
\end{remark}

\section{Perturbed hierarchies and fluctuation limits}\label{S:fluctuation}

We very briefly discuss fluctuation results for reaction diffusion systems. Precise results and details can be found in \cite{pD88b, zqCwtF14c, zqCwtF14d}.

The fluctuation of the empirical measure $\X^N$  in \eqref{Def:EmpiricalRD} at time $t$ is defined by
$$\Y^N_t(\phi) := \sqrt{N}\,(\<\X^N_t,\phi\>-\E\<\X^N_t,\phi\>).$$
Even in our simple setting in Section \ref{S:model}, it is nontrivial to obtain satisfactory answers to the following natural questions:
\begin{enumerate}
	\item[(1)]   What is the state space for $\Y^N_t$? This space should posses a topology which allows us to make sense of convergence of $\Y^N$. 
	\item[(2)] If it does converge, how to prove convergence and what can we say about the limit $\Y$?
\end{enumerate}

For fluctuation results, the case $R(u)=-u^2$ is treated in \cite{pD88b}, and more general cases in
\cite{pK86, pK88, BDP92}. Roughly speaking, the fluctuation limit $\Y$   solves the following {\it stochastic partial differential equation} in a distributional Hilbert space:
\begin{equation*}
d\Y_t = \Big(\frac{1}{2}\Delta\Y_t+R'(u(t))\Y_t\Big)\,dt + dM_t,
\end{equation*}
where $u(t,x)$ solves equation \eqref{E:RDeqt}, $R'(u)$ is the derivative of $R(u)$ (e.g. $-2u$ when $R(u)=-u^2$) and is viewed as a multiplicative operator, $M$ is a Gaussian martingale with independent increment and covariance structure
\begin{equation}\label{CovM}
\E[(M_t(\phi))^2]=\int_0^t\<|\nabla\phi|^2,\,u(s)\>+\<\phi^2,\,|R(u(s))|\>\,ds.
\end{equation} 
Here $\<\cdot\,,\,\cdot\>$ is the $L^2$ inner product in the spatial variable and $|R(u)|$ is the polynomial obtained by putting an absolute sign to each coefficient in $R(u)$. Observe that \eqref{CovM} is the formal limit of \eqref{E:Mtg_phi_quad} after taking expectation.
%These results are restricted either to polynomials of low order or to lower dimensions.  
%The functional analytic framework for fluctuation results in  domains is developed in \cite{zqCwtF14c}. Boundary conditions are encoded in the function spaces in consideration.

In \cite{pD88b}, a key step in establishing a fluctuation result for  the case $R(u)=-u^2$ is to compute  the  second order  approximation of the correlation function
$F^{(k)}$. That is, one find out an expression for $G^{N,(k)}$ such that
\begin{equation}\label{E:2nd}
F^{(k)}_t(\vec{x})= \prod_{i=1}^k u(t,x_i)+\frac{G^{N,(k)}_t(\vec{x})}{N} + \frac{o(N)}{N}.
\end{equation}

The key idea is to regard the terms $\frac{QF^{(k)}_s}{N}$ in \eqref{E:BBGKY_F} as "small errors" and {\it introduce two approximating hierarchies}
\begin{align*}
A^{(k)}_t&=P^{(k)}_t F^{(k)}_0 - \int_0^t P^{(k)}_{t-s}\big( RA^{(k+1)}_s \big)\;ds \quad\text{and}\\
B^{(k)}_t&=P^{(k)}_t F^{(k)}_0 - \int_0^t P^{(k)}_{t-s}\left( RB^{(k+1)}_s + \frac{QA^{(k)}_s}{N} \right)\;ds.
\end{align*}
The remarkable point is that these two hierarchies have explicit product form solutions when $F^{(k)}_0$ has product form. In fact, $A^{(k)}_t(\vec{x}) =\prod_{i=1}^k u_N(t,x_i)$ where $u_N$ is uniformly close to $u$,  $B^{(k)}_t=A^{(k)}_t + \frac{G^{N,(k)}_t}{N}$ for an explicit function $G^{N,(k)}_t$ and $N\,(F^{(k)}_t-B^{(k)}_t)\to 0$ uniformly in $(t,\vec{x})$ in a compact set as $N\to\infty$. Hence we obtain \eqref{E:2nd}. See \cite{pD88b} for details.

\begin{remark}\label{Rk:domain}\rm (Extensions)
	Here we focus on the case 
	$\bar{D}=[0,1]$, but the theorem and the proof can be generalized to any bounded Lipschitz domain. The essential point is that the domain of the Feller generator of the RBM on $D$, denoted by $Dom(\frac{1}{2}\Delta)$, is dense in $C(\bar{D})$.
	One can also prove functional LLN for reaction diffusion systems without going through the BBGKY. See, for instance, the perturbation method in \cite{pD88a} or trick of interchanging limits in \cite{zqCwtF14b}.
	It is well-known that for exchangeable systems, propagation of chaos is equivalent to LLN. The precise statement can be found in \cite{Szn91}. The  functional LLN result works for non-exchangeable systems as well. 
		%Theorem \ref{T:POC_RD} does not directly imply Theorem \ref{T:LLN_RD} since the former is not a pathwise result; one also needs to prove tightness in the path space to fill this gap. 
\end{remark}


\begin{thebibliography}{99} %\footnotesize
  \bibitem{BDP92}
  C. Boldrighini, A. De Masi and A. Pellegrinotti.
  Nonequilibrium fluctuations in particle systems modelling reaction-diffusion equations.
  \emph{Stochastic Processes. Appl.} \textbf{42} (1992), 1-30.

  \bibitem{BFW}
A. Budhiraja, W.-T. Fan and R. Wu.
Large deviations for Brownian particle systems with killing.  {\it  Journal of Theoretical Probability}.  (2016) 1-40.


  \bibitem{zqCwtF14b}
  Z.-Q. Chen and W.-T. Fan.
  Systems of interacting diffusions with partial annihilations through membranes.  {\it Annals of Probability}. \textbf{45} (2017) 100-146.

  \bibitem{zqCwtF13a}
  Z.-Q. Chen and W.-T. Fan.
  Hydrodynamic limits and propagation of chaos for interacting random walks in domains.  {\it Annals of Applied Probability}.  \textbf{27}(3) (2017) 1299-1371.

  \bibitem{zqCwtF14c}
  Z.-Q. Chen and W.-T. Fan.
  Functional central limit theorem for Brownian particles in domains with Robin boundary condition. 
  \emph{Journal of functional analysis.}  \textbf{269}(12) (2015), 3765-3811.
  
  \bibitem{zqCwtF14d}
  Z.-Q. Chen and W.-T. Fan.
  Fluctuation limit for systems of interacting diffusions with partial annihilations through membranes 
  \emph{Journal of Statistical Physics.}  \textbf{164} (2016), 890—936.
 

	\bibitem{MP91}
	De Masi, A. and Presutti, E.  
	{\it Mathematical methods for hydrodynamic limits.} 
	Lecture Notes in Mathematics. 1991.

  \bibitem{pD88a}
    P. Dittrich.
    A stochastic model of a chemical reaction with diffusion.
    \emph{Probab. Theory Relat. Fields.} \textbf{79} (1988), 115-128.
  \bibitem{pD88b}
    P. Dittrich.
    A stochastic partical system: Fluctuations around a nonlinear reaction-diffusion equation.
    \emph{Stochastic Processes. Appl.} \textbf{30} (1988), 149-164.

  \bibitem{lEbShtY07}
  L. Erd\"os, B. Schlein and H. T. Yau.
  Derivation of the cubic non-linear Schr\"odinger equation from quantum dynamics of many-body systems.
  \emph{Inventiones Mathematicae.} \textbf{167}(3) (2007), 515-614.

  \bibitem{EK86}
  S.N. Ethier, S.N. and Kurtz, T.G.
  {\it Markov processes. Characterization and Convergence.}
  Wiley, New York, 1986. MR0838085.



  \bibitem{pK86}
  P. Kotelenez.
  Law of large numbers and central limit theorem for linear chemical reactions with diffusion.
  \emph{Ann. Probab.} \textbf{14} (1986), 173-193
  
  \bibitem{pK88}
  P. Kotelenez.
  High density limit theorems for nonlinear chemical reactions with diffusion.
  \emph{Probab. Theory Relat. Fields.} \textbf{78} (1988), 11-37.


  \bibitem{Sznitman87}
A. S. Sznitman
Propagation of chaos for a system of annihilating Brownian spheres.
\emph{Commm. Pure Appl. Math.} \textbf{6} (1987), 663-690.
 
\bibitem{Szn91}
A.-S. Sznitman. Topics in propagation of chaos. 
\textit{Lecture Notes in Mathematics} Vol 1464, 165--251, 1991.


\end{thebibliography}
\end{document}